\documentclass{amsart}
\usepackage{amsmath,amssymb}
\usepackage[matrix,arrow,cmtip]{xy}
\newcommand{\PP}{\mathbb{P}}
\newcommand{\ZZ}{\mathbb{Z}}
\newcommand{\QQ}{\mathbb{Q}}

\newcommand{\Disc}{\mathrm{Disc}}

\newcommand{\tors}{\mathrm{tors}}
\newcommand{\rk}{\mathrm{rk}}

\newcommand{\fp}{\mathfrak{p}}

\newlength{\underscorelength}
\renewcommand{\phi}{\varphi}
\newtheorem{theorem}{Theorem}[section]
\newtheorem{lemma}[theorem]{Lemma}

\theoremstyle{definition}

\begin{document}
\title{The primitive solutions to $x^3+y^9=z^2$}
\parskip=0.1pt
\author{Nils Bruin}
\thanks{The author is partially supported by NSERC and the research
described in this paper was done partially at the University of Sydney}
\address{Department of Mathematics, Simon Fraser University,
Burnaby BC, Canada V5A 1S6}
\email{bruin@member.ams.org}


\date{October 31, 2003}

\begin{abstract}

We determine the rational integers $x,y,z$ such that $x^3+y^9=z^2$ and
$\gcd(x,y,z)=1$. First we determine a finite set of curves of genus $10$ such
that any primitive solution to $x^3+y^9=z^2$ corresponds to a rational point on
one of those curves. We observe that each of these genus $10$ curves covers an
elliptic curve over some extension of $\QQ$. We use this cover to apply a
Chabauty-like method to an embedding of the curve in the Weil restriction of the
elliptic curve. This enables us to find all rational points and therefore deduce
the primitive solutions to the original equation.

\end{abstract}
\maketitle
\section{Introduction}

In this article, we consider a special instance of the equation
$Ax^r+By^s=Cz^t$. One of the most important results for this equation is a
theorem  by Darmon and Granville (\cite{DG:genferm}), which for
fixed, nonzero $A,B,C$ and fixed $r,s,t$ satisfying $1/r+1/s+1/t<1$, relates
primitive solutions $x,y,z$ (solutions with $x,y,z\in\ZZ$ and $\gcd(x,y,z)=1$)
to rational points on finitely many algebraic
curves of general type. It follows by Faltings' theorem that only finitely many
primitive solutions exist. Unfortunately, their proof does not provide a recipe
to produce the relevant curves.

As Tijdeman shows in \cite{tij:banff}, the ABC-conjecture implies that if one
leaves $A,B,C$ fixed, but allows $r,s,t$ to vary subject to $1/r+1/s+1/t<1$,
then the number for primitive solutions to $Ax^r+By^s=Cz^t$ is still finite.

The case $A=B=C=1$ has received most attention. An extensive computer search
performed by Beukers and Zagier (see \cite{beu:dioph}) showed there are some
surprisingly large primitive solutions to $x^r+y^s=z^t$ for
\begin{equation}\label{eq:nontrivexp}
\{r,s,t\}\in\{\{2,3,7\},\{2,3,8\},\{2,3,9\},\{2,4,5\}\}.
\end{equation}
We call a solution \emph{trivial} if one of $x,y,z$ is equal to $\pm 1$.
Note that this includes $2^3+1^r=3^2$ and, for primitive solutions, also any
solution satisfying $xyz=0$. With this definition, the only equations
$x^r+y^s=z^t$ with $1/r+1/s+1/t<1$ for which non-trivial primitive solutions are
known, are the ones satisfying (\ref{eq:nontrivexp}).

It was noted by Tijdeman and Zagier that the known non-trivial primitive
solutions of $x^r+y^s=z^t$ all have $\min(r,s,t)\leq 2$. They conjectured that 
for $r,s,t\geq 3$ only trivial primitive solutions exist. The resolution of this
conjecture has even attracted a monetary prize (see \cite{mauldin:beal}).

In this paper, we work in a different direction. Because the equations
satisfying (\ref{eq:nontrivexp}) have known non-trivial primitive solutions,
the question whether they have any others appears natural. In
his thesis (\cite{bruin:phdthesis}, \cite{bruin:crelle} and \cite{bruin:283}),
the author resolved the cases
$\{r,s,t\}\in\{\{2,3,8\},\{2,4,5\}\}$.
Recently Poonen, Schaefer and Stoll have made excellent progress on
$x^2+y^3=z^7$. In this paper we deal with the one remaining case:

\begin{theorem}\label{thrm:main}
The primitive solutions to the equation $x^3+y^9=z^2$ are
$$(x,y,z)\in\{(1,1,0),(0,1,\pm1),(1,0,\pm1),(2,1,\pm3),(-7,2,\pm13)\}.$$
\end{theorem}

\section{Parametrising curves for $x^3+y^9=z^2$}

We note that the perfect $9$th powers form a subset of the perfect cubes, so any
primitive solution $(x,y,z)$ to $x^3+y^9=z^2$ gives rise to a primitive solution
$(x,v,z)$ of $x^3+v^3=z^2$ by putting $v=y^3$. The solutions of the latter
equation can be classified by the following theorem.

\begin{theorem}[Mordell, {\cite[Chapter 25]{mordell:diopheq}}]
\label{thrm:233}
If $x,v,z\in\ZZ$ with $x^3+v^3=z^2$ and $\gcd(x,v,z)=1$, then
there are $s,t\in\ZZ[\frac{1}{2},\frac{1}{3}]$ and
$(s,t)\ne(0,0)\mod p$ for any prime $p\nmid 6$, such one of the following holds.
$$\begin{array}{rcl}
x\mbox{ or }v&=&s^4+6 s^2 t^2-3 t^4\\
v\mbox{ or }x&=&-s^4+6 s^2 t^2+3 t^4\\
z&=&\pm 6 s t (s^4+3 t^4)\\
\end{array}$$
$$\begin{array}{rcl}
x\mbox{ or }v&=&\frac{1}{4} (s^4+6 s^2 t^2-3 t^4)\\
v\mbox{ or }x&=&\frac{1}{4} (-s^4+6 s^2 t^2+3 t^4)\\
z&=&\pm \frac{3}{4} s t (s^4+3 t^4)\\
\end{array}$$
$$\begin{array}{rcl}
x\mbox{ or }v&=&s (s^3+8 t^3)\\
v\mbox{ or }x&=&4 t (t^3-s^3)\\
z&=&\pm s^6-20 s^3 t^3-8 t^6\\
\end{array}$$
\end{theorem}
It follows that in order to find all primitive solutions to $x^3+y^9=z^2$, it
suffices to find all solutions $y,s,t$ with $y\in\ZZ$ and
$s,t\in\ZZ[\frac{1}{2},\frac{1}{3}]$ with $s,t$ not both divisible by any prime
$p$ not dividing 6, to the following equations:
$$\begin{array}{rrcl}
1:&y^3&=&s^4+6 s^2 t^2-3 t^4\\
2:&y^3&=&-s^4+6 s^2 t^2+3 t^4\\
3:&y^3&=&\frac{1}{4} (s^4+6 s^2 t^2-3 t^4)\\
4:&y^3&=&\frac{1}{4} (-s^4+6 s^2 t^2+3 t^4)\\
5:&y^3&=&s (s^3+8 t^3)\\
6:&y^3&=&4 t (t^3-s^3).
\end{array}$$
This leads us to a generalisation of the concept of \emph{primitive solution}.
Let $S=\{p_1,\ldots,p_s\}$ be a finite set of primes. We write
$$\ZZ_S=\ZZ[1/p_1,\ldots,1/p_s].$$
A tuple $(x_1,\ldots,x_n)\in\ZZ_S$ is called \emph{$S$-primitive} if
the ideal $(x_1,\ldots,x_n)\ZZ_S$ equals $\ZZ_S$. Equivalently, it means that
no prime outside $S$ divides all $x_i$. In order to determine the primitive
solutions of $x^3+y^9=z^2$, it suffices to determine the \emph{$S$-primitive}
solutions to the equations above.

Obviously, if $(s,t,y)$ is a solution to one of these equations, then so is an
entire class of weighted homogeneously equivalent solutions of the form
$(\lambda^3s,\lambda^3t,\lambda^4y)$. Furthermore, a solution $(s,t,y)$ to one
of the equations 1 or 2 gives rise to a solution
$(\frac{1}{2}s,\frac{1}{2}t,\frac{1}{4}y)$ of equations $3$ or $4$
respectively. Note that the quantity $(s:t)$ is invariant under either
transformation. Since a solution can be easily reconstructed (up to
equivalence), it is sufficient to determine the values of $s/t$ that can occur
for $\{2,3\}$-primitive solutions of the equations 1, 2, 5 and 6.

Following \cite{DG:genferm}, the equivalence classes under homogeneous
equivalence of $\{2,3\}$-primitive solutions to each of these equations
correspond to the rational points on a finite collection of algebraic curves.
We construct these curves explicitly (see also \cite{bruin:crelle}). We
introduce some notation.

Let $K$ be a number field and let $S$ be a finite set of rational primes. For a prime
$\fp$ of $K$ we write $\fp\nmid S$ if $v_\fp(q)= 0$ for all $q\in S$. Following
\cite{sil:AEC1}, for any rational prime $p$ we define the following subgroup
of $K^*/K^{*p}$:
$$K(p,S)=\{a\in K^*:
   v_\fp(a)\equiv 0 \pmod{r} \mbox{ for all primes }\fp\nmid S\}/K^{*p}.$$
This is a finite group. We will identify it with a set of representatives in
$K^*$.

Let $f\in K[x]$ be a square-free polynomial and let $A:=K[x]/(f)$.
The algebra $A$ is isomorphic to a direct sum of
number fields $K_1,\ldots,K_r$ and $A^*=K_1^*\times\cdots\times K_r^*$. We
generalise the notation above by defining
$$A(p,S):=K_1(p,S)\times\cdots\times K_r(p,S)\subset A^*/A^{*p}.$$
We identify the elements of this finite group with a set of representatives in
$A^*$.

\begin{lemma}\label{lemma:genparm} Let $f(s,t)\in\ZZ[s,t]$ be a square-free homogeneous form of
degree $4$. Let $S=\{p \mbox{ prime}:p \mid 3\Disc(f)\}$. Solutions of
$y^3=f(s,t)$ with $y,s,t\in\QQ$, $s,t$ integral outside $S$ and
$(s,t)\mod p\neq (0,0)$ for any $p\notin S$ correspond, up to weighted
projective equivalence, to rational points on
finitely many explicitly constructible smooth projective curves of the form
$$\begin{array}{rcl}
Q_{2,\delta}(y_0,y_1,y_2,y_3)&=&0\\
Q_{3,\delta}(y_0,y_1,y_2,y_3)&=&0,
\end{array}$$
indexed by $\delta\in A(3,S)$ for some explicit semisimple algebra $A$.
For each of these curves, the ratio $s/t$ can be explicitly expressed as a
function
$$\frac{s}{t}=-
   \frac{Q_{0,\delta}(y_0,y_1,y_2,y_3)}{Q_{1,\delta}(y_0,y_1,y_2,y_3)}.$$
\end{lemma}
\begin{proof}
Using an $\mathrm{SL}_2(\ZZ)$ transformation on $s,t$ if necessary, we can assume that
$f$ is of degree $4$ in $s$. 
We form the algebra $A=\QQ[\theta]=\QQ[x]/f(x,1)$. Let $C$ be the leading
coefficient of $f(x,1)$.
Then we have $f=C N_{A[s,t]/\QQ[s,t]}(s-\theta t)$, so having a solution would
amount to the existence of $y_0,y_1,y_2,y_3\in\QQ$ and $\delta\in A^*$ such that
\begin{eqnarray*}
(s-\theta t)&=&\delta(y_0+\theta y_1+\theta^2 y_2+\theta^3 y_3)^3\\
y&=&\sqrt[3]{C\,N_{A/\QQ}(\delta)}\,
N_{A[s,t]/\QQ[s,t]}(y_0+\theta y_1+\theta^2 y_2+\theta^3 y_3).
\end{eqnarray*}
It follows immediately that $C\,N_{A/\QQ}(\delta)$ should be a perfect cube and
that the value of $\delta$ is only relevant modulo cubes.

The fact that $s,t$ are coprime implies that $\delta$ can be chosen from a
finite set of representatives.  From, for instance,
\cite[Lemma~2.2.1]{bruin:phdthesis} it follows that if $s,t$ are integral and
coprime outside of $S$, then $\delta\in A(3,S)$.

Given $\delta\in A(3,S)$ with $N(\delta)\in\QQ^{*3}$, there are unique cubic
forms $Q_{0,\delta},\ldots,Q_{3,\delta}\in\QQ[y_0,\cdots,y_3]$ such that
$$\delta(y_0+\theta y_1+\theta^2 y_2+\theta^3 y_3)^3=
Q_{0,\delta}+Q_{1,\delta}\theta+Q_{2,\delta}\theta^2+Q_{3,\delta}\theta^3$$
Solving $s,t$ from the equations above, we obtain
$$\begin{array}{rcl}
s&=&Q_{0,\delta}(y_0,y_1,y_2,y_3)\\
t&=&-Q_{1,\delta}(y_0,y_1,y_2,y_3)\\
0&=&Q_{2,\delta}(y_0,y_1,y_2,y_3)\\
0&=&Q_{3,\delta}(y_0,y_1,y_2,y_3)\\
\end{array}$$
The latter two equations define a smooth projective curve of genus $10$
(see \cite{bruin:phdthesis}) and the first two equations show that $s/t$ is a
rational function on that curve.
\end{proof}

\section{Rational points on the parametrising curves}
We now proceed to apply this procedure to each of the equations mentioned above
and then determine the finite set of values that $s/t$ attains for the
$\{2,3\}$-primitive solutions. Lemmas~\ref{lemma:eq5} and \ref{lemma:eq12} both
require standard but elaborate and tedious computations. These are easily
executed by a computer algebra system, but are too bulky to reproduce completely
on paper. We only give the basic data to check these computations. For the
interested reader, we have made available a complete electronic transcript
\cite{bruin:comp239} of the computations in MAGMA \cite{magma}.

\begin{lemma}\label{lemma:eq5}
The
$\{2,3\}$-primitive solutions of the equation
$$y^3=s (s^3+8 t^3)$$
have
$$\frac{s}{t}\in\{-2,0,1,2,4,\infty\}.$$
\end{lemma}
\begin{proof}
We use the construction in the proof of Lemma~\ref{lemma:genparm}. In our case,
$A=\QQ[\theta]$, where $\theta(\theta+2)(\theta^2-2\theta+4)=0$ and $S=\{2,3\}$.
The subgroup of $A(3,S)$ of elements of cubic norm is spanned by
$$\begin{array}{l}
    -1/8\theta^3 - 1/4\theta^2 - \theta + 1,\\
    -1/24\theta^3 - 1/6\theta^2 - 2/3\theta + 1,\\
    1/24\theta^3 + 1/6\theta^2 + 1/6\theta + 1,\\
    5/24\theta^3 - 5/12\theta^2 - 2/3\theta + 3,\\
    1/8\theta^3 - 1/4\theta^2 - 1/2\theta + 2.
\end{array}$$
For each of the $3^5$ possible values of $\delta$ in that group, we can write
down the curve $Q_{2,\delta}=Q_{3,\delta}=0$ as outlined in the proof of
Lemma~\ref{lemma:genparm} and test the curve for local solubility over $\QQ_3$. Only
$22$ values for $\delta$ pass this test.

From the factorisation of $f$ it follows that there are two homomorphisms
$m_1,m_2:A\to\QQ$, defined by $m_1(\theta)=0$ and $m_2(\theta)=-2$. For a fixed
value of $\delta$, the curve $Q_{2,\delta}=Q_{3,\delta}=0$ covers two curves of
genus $1$:
$$\begin{array}{ccccc}
E_{1,\delta}&:& \frac{N(\delta)}{m_1(\delta)}u_1^3 &=&s^3-8t^3\\
E_{2,\delta}&:& \frac{N(\delta)}{m_2(\delta)}u_2^3 &=&s(s^2-2st+4t^2)
\end{array}$$
\begin{table}
\begin{small}
$$\begin{array}{c|c|c|c|c}
\delta&N(\delta)/m_1(\delta)&N(\delta)/m_2(\delta)&\rk E_{1,\delta}(\QQ)&
\rk E_{2,\delta}(\QQ)\\
\hline
1&1&1&0&1\\
-\frac{3}{8}\theta^3+\frac{1}{2}\theta^2-\frac{3}{2}\theta+1&1&3&0&0\\
\frac{1}{8}\theta^3-3\theta+1&1&36&0&1\\
\frac{1}{8}\theta^3-2\theta+1&1&2&0&0\\
\frac{1}{24}\theta^3+\frac{1}{6}\theta^2+\frac{1}{6}\theta+1&1&1&0&1\\
-\frac{3}{8}\theta^3+\frac{3}{4}\theta^2-\theta+1&1&3&0&0\\
-\frac{7}{8}\theta^3+\theta^2-3\theta+1&1&12&0&0\\
\frac{1}{8}\theta^3-6\theta+1&1&18&0&1\\
-\frac{1}{8}\theta^3+\frac{1}{2}\theta^2+\frac{1}{2}\theta+1&1&9&0&0\\
-\frac{1}{24}\theta^3+\frac{1}{3}\theta^2+\frac{1}{3}\theta+1&1&4&0&0\\
-\frac{15}{8}\theta^3+2\theta^2-6\theta+1&1&6&0&0\\
-\frac{1}{12}\theta^3+\frac{2}{3}\theta^2-\frac{4}{3}\theta+3&9&3&1&0\\
-\frac{23}{24}\theta^3+\frac{8}{3}\theta^2-\frac{22}{3}\theta+3&9&6&1&0\\
\frac{1}{24}\theta^3+\frac{11}{12}\theta^2-\frac{4}{3}\theta+3&9&3&1&0\\
\frac{5}{8}\theta^3+\frac{1}{4}\theta^2-2\theta+9&3&3&0&0\\
-\frac{3}{8}\theta^3+\frac{3}{4}\theta^2-\frac{1}{2}\theta+2&4&3&0&0\\
\frac{7}{24}\theta^3+\frac{5}{12}\theta^2-\frac{11}{6}\theta+6&36&3&0&0\\
\frac{13}{8}\theta^3-\frac{1}{4}\theta^2-\frac{5}{2}\theta+18&12&3&1&0\\
\frac{5}{8}\theta^3+\frac{3}{4}\theta^2-\frac{7}{2}\theta+4&2&3&0&0\\
\frac{7}{24}\theta^3-\frac{7}{12}\theta^2-\frac{5}{6}\theta+4&2&1&0&1\\
\frac{23}{24}\theta^3+\frac{1}{12}\theta^2-\frac{13}{6}\theta+12&18&3&0&0\\
\frac{29}{8}\theta^3-\frac{5}{4}\theta^2-\frac{7}{2}\theta+36&6&3&1&0\\
\end{array}$$
\end{small}
\caption{Ranks of $E_{i,\delta}(\QQ)$ related to $y^3=s(s^3+8t^3)$}
\label{tbl:1}
\end{table}
Note that $(s:t:u_1)=(2:1:0)\in E_{1,\delta}(\QQ)$ and $(s,t,u_2)=(0:1:0)\in
E_{2,\delta}(\QQ)$, so both curves are isomorphic to elliptic curves and hence
their rational points have the structure of a finitely generated group.
For the $22$ values of $\delta$, we obtain the results in Table~\ref{tbl:1}.
For each of the values of $\delta$ at least one of
$E_{1,\delta}(\QQ)$ and $E_{2,\delta}(\QQ)$ is of rank $0$. Hence, any
$\{2,3\}$-primitive solution must have $s/t$ corresponding to a torsion point on
one of $E_{1,\delta}(\QQ)$ and $E_{2,\delta}(\QQ)$.

The group $E_{1,\delta}(\QQ)^\tors$ is $\ZZ/3\ZZ$ for
$N(\delta)/m_1(\delta)=1$ and $\ZZ/2\ZZ$ for $N(\delta)/m_1(\delta)=2$.
For other values there is no torsion. The non-trivial $3$-torsion points are
$(s:t:u_1)=(1:0:1),(0:1:2)$. The non-trivial $2$-torsion point is
$(s:t:u_1)=(2:1:8)$.

The group $E_{2,\delta}(\QQ)^\tors$ is $\ZZ/3\ZZ$ for
$N(\delta)/m_1(\delta)=3$ and $\ZZ/2\ZZ$ for $N(\delta)/m_1(\delta)=6$. For
other values there is no torsion. The non-trivial $3$-torsion points are
$(s:t:u_2)=(1:1:1),(2:-1:8)$. The non-trivial $2$-torsion point is
$(s:t:u_2)=(4:1:2)$. These points give rise to the set of values stated in the
lemma.
\end{proof}

\begin{lemma}\label{lemma:eq6}
 The
$\{2,3\}$-primitive solutions of the equation
$$y^3=4t (t^3-s^3)$$
have
$$\frac{s}{t}\in\{-2,-1,-\frac{1}{2},0,1,\infty\}.$$
\end{lemma}
\begin{proof}
Note that the map $(s,t,y)\mapsto(-t/2,s/4,y/4)$ is a bijection from
the $\{2,3\}$-primitive solutions of $y^2=s(s^3+8t^3)$ to those of
$y^3=4t(t^3-s^3)$. Lemma~\ref{lemma:eq5} together with the induced map
$s/t\mapsto -2t/s$ proves the statement.
\end{proof}

\begin{lemma}\label{lemma:eq12}
The $\{2,3\}$-primitive solutions of
$$y^3=s^4+6 s^2 t^2-3 t^4$$
have
$$\frac{s}{t}\in\{-1,0,1,\infty\}$$
and the $\{2,3\}$-primitive solutions of
$$y^3=-s^4+6 s^2 t^2+3 t^4$$
have
$$\frac{s}{t}\in\{-3-1,0,1,3,\infty\}.$$
\end{lemma}
\begin{proof}
Again we follow the procedure outlined in Lemma~\ref{lemma:genparm}. For both
equations, the algebra $A$ is isomorphic to the number field
$K(\alpha):=\QQ[x]/(x^4-2x^3-2x+1)$. In each case, we are left with $4$ values of $\delta$
for which the curve $C_\delta:Q_{\delta,2}=Q_{\delta,3}=0$ has points over
$\QQ_3$. Each
of those curves actually has a rational point. Over $K$, it covers the genus
$1$ curve
$$E_\delta:\frac{N(\delta)}{\delta}u^3=\frac{N(s-\theta t)}{s-\theta t}$$
The existence of a $\QQ$-rational point on $C_\delta$ implies that there is a
$K$-rational point on
$E_\delta$, which makes $E_\delta$ isomorphic to an elliptic curve. We have the
following diagram.
$$\xymatrix@R-2em{
C_{/\QQ} \ar[rd]^{\pi} \ar[dd]_{\frac{s}{t}}\\
&E_{/K}\ar[dl]^{\frac{s}{t}}\\
\PP^1_{/\QQ}}$$
It follows immediately that
$$\frac{s}{t}(C(\QQ))\subset\PP^1(\QQ)\cap\frac{s}{t}(E(K)).$$
A method for determining the intersection on the right-hand side is described
in \cite{bruin:crelle}. The method is an adaptation of Chabauty's method
(\cite{chab:ratpoints}) and applies if $\rk E(K) < [K:\QQ]$. The method requires
generators of a subgroup of $E(K)$ of maximal rank. Table~\ref{tbl:mwgrp} lists such
generators for the six elliptic curves that are encountered. A $2$-descent shows
that these points generate a subgroup of finite index prime to $2$ and a simple
check for $3$-divisibility of all relevant linear combinations of these points
that the index is also prime to $3$. We will need this later on.
\begin{table}
\begin{small}
$$\begin{array}{l|l|l}
i&E_i:y^2=\ldots&\mbox{independent points in }E_i(K)\\
\hline
1&x^3 + 2\alpha^3 - 4\alpha^2 - 2\alpha - 6&
        g_1=(2,\alpha^3 - 2\alpha^2 - \alpha)\\
      &&g_2=(2\alpha^2 + 2\alpha + 2,-3\alpha^3 - 4\alpha^2 - 5\alpha - 2)\\
\hline
2&x^3 - 2\alpha^3 + 4\alpha^2 + 2\alpha - 2&
        (2\alpha,\alpha^3 - \alpha)\\
      &&(1,\alpha^2 - 2\alpha)\\
\hline
3&x^3 + 2\alpha^3 - 6\alpha^2 + 2\alpha&
        (2,\alpha^2 - 3)\\
\hline
4&x^3 + 2\alpha^3 - 2\alpha^2 - 6\alpha&
        (2,2\alpha^3 - 3\alpha^2 - 2\alpha - 1)\\
\hline
5&x^3 - 30\alpha^3 - 22\alpha^2 - 30\alpha&
        (2\alpha^2 + 2\alpha + 2,6\alpha^3 - \alpha^2 + 2\alpha - 5)\\
\hline
6&x^3 + 14\alpha^3 - 6\alpha^2 - 74\alpha + 32&
        (-2\alpha^3 + 6\alpha^2 - 4\alpha + 2
                6\alpha^3 - 13\alpha^2 - 4\alpha + 5)
\end{array}$$
\end{small}
\caption{Independent points on $E_i$ over $K=\QQ(\alpha)$}
\label{tbl:mwgrp}
\end{table}
For the first equation, using $\theta=\alpha^2-2\alpha$, we find that
$C_\delta(\QQ_3)\neq \emptyset$ for the $4$ values of $\delta$ listed in
Table~\ref{tbl:deltas}. We
also list the value of $\frac{s}{t}$ at a point $p_0\in C_\delta(\QQ)$ and we
list $i$ such that $E_\delta$ is isomorphic to $E_i$ in Table~\ref{tbl:mwgrp}.
\begin{table}
\begin{small}
\begin{tabular}{ccc}
$\begin{array}{c|c|l}
\delta&\frac{s}{t}(p_0)&i\\
\hline
1&\infty&1\\
\alpha^3 - 2\alpha^2 - \alpha - 3&0&2\\
 -3\alpha^3 + 5\alpha^2 + 2\alpha + 7&1&3\\
7\alpha^3 - 11\alpha^2 - 5\alpha - 16&-1&4
\end{array}$&\hspace{1em}&
$\begin{array}{c|c|l}
\delta&\frac{s}{t}(p_0)&i\\
\hline
    1&\infty&2\\
    -4\alpha^3 + 8\alpha^2 + 4\alpha + 11&0&1\\
    \alpha^2 + \alpha + 1&1&6\\
    \alpha^3 + \alpha^2 + \alpha&-1&5
\end{array}$\\
$\theta=\alpha^2-2\alpha$&&$\theta=\alpha^3 - \alpha^2 - \alpha - 2$
\end{tabular}
\end{small}
\caption{$E_\delta$, $\frac{s}{t}(p_0)$ and isomorphic $E_i$}
\label{tbl:deltas}
\end{table}
For the second equation, using $\theta=\alpha^3 - \alpha^2 - \alpha - 2$, we find
a similar set of values.
Note that all the eight smooth plane cubics $E_\delta$ have a $K$-rational point
of inflection. Thus, the isomorphism between $E_\delta$ and $E_i$ can be realised
by a linear change of variables. As an example, we give some data for
$\delta=-4\alpha^3 + 8\alpha^2 + 4\alpha + 11$. The point $(u:s:t)=
(0:-\alpha^3 + \alpha^2 + \alpha + 2 : 1)$ is a point of inflection. If we map
this point to the origin on $E_1$, then the function $\frac{s}{t}$ on $E_1$ can
be expressed as
$$\frac{s}{t}=\frac{(-\alpha^3 + \alpha^2 + \alpha + 2)y+
(-\alpha^3 + 2\alpha^2 - \alpha + 2)}{y+(3\alpha^3 - 6\alpha^2 - 3\alpha - 8)},$$
where $y$ is a coordinate of the Weierstrass model of $E_1$ as mentioned in
Table~\ref{tbl:mwgrp}.

We find $\frac{s}{t}(\{g_1,-g_2,-g_1+g_2\})=\{0,-3,3\}$. Via a Chabauty-like
argument using the primes of $K$ over $11$ (see \cite{bruin:crelle}), we find
that these are all points in $E_1(K)$ that have values in $\QQ$ under
$\frac{s}{t}$.

For the other curves we can use the same method. Either the primes over $11$ or
over $31$ yield sufficient information to conclude that there is only one point
with a $\QQ$-rational value of $\frac{s}{t}$. This is the value listed in the
tables above.
\end{proof}

\section{Proof of Theorem~\ref{thrm:main}}
Given Lemmas~\ref{lemma:eq5}, \ref{lemma:eq6} and \ref{lemma:eq12}, the proof of
Theorem~\ref{thrm:main} is reduced to filling in pairs of $s,t$ giving rise to
the values $\frac{s}{t}$ listed in those lemmas in the corresponding
parametrisations listed in Theorem~\ref{thrm:233}. We then check
if the resulting solution to $x^3+v^3=z^2$ is weighted homogeneously
equivalent with a primitive solution to $x^3+y^9=z^2$.

Since the Lemmas~\ref{lemma:eq5} and \ref{lemma:eq6} give
values of $\frac{s}{t}$ for wich the first or the second polynomial of the third
parametrization in Theorem~\ref{thrm:233} may be a cube, we try all of the
values and see what we get.
$$\begin{array}{rr|c}
s&t&\mbox{solution}\\
\hline
-2&1&0^3+36^3=216^2\\
0&1&0^3+4^3=(-8)^2\\
1&1&9^3+0^3=(-27)^2\\
1&0&1^3+0^3=1^3\\
2&1&32^3-28^3=(-104)^2\\
4&1&288^3-252^3=2808^2\\
-1&1&-7^3+8^3=13^2\\
-1&2&-63^3+72^3=(-351)^2
\end{array}$$
It is easily checked that all solutions are equivalent to one occurring in the
list in the theorem. Note that the first four entries correspond to values of
$\frac{s}{t}$ that occur in both Lemma~\ref{lemma:eq5} and \ref{lemma:eq6}.
Correspondingly, both cubes are ninth powers, up to $\{2,3\}$-units. The latter
four values occur in only one of the lemmas and correspondingly have primes
different from $2,3$ dividing one of the cubes to a power not divisible by 9.

For the values of $\frac{s}{t}$ listed in
Theorem~\ref{lemma:eq12} we proceed similarly.
$$\begin{array}{rr|c}
s&t&\mbox{solution}\\
\hline
1&0&1^3-1^3=0^2\\
0&1&-3^3+3^3=0^2\\
\pm1&1&4^3+8^3=(\pm24)^2\\
\pm3&1&132^3-24^3=(\pm1512)^2
\end{array}$$
The last two solutions are interesting. The penultimate solution is obviously
equivalent to $1^3+2^3=3^2$ and, since $1$ is indeed a cube, this yields a
solution. The last one is equivalent to $33^3-6^3=189^2$. Indeed, this is
equivalent to the solution
$(2^2\cdot3^3\cdot11)^3-(2\cdot3)^9=(2^3\cdot3^6\cdot7)^2$, which is a
$\{2,3\}$-primitive solution. It is clear, however, that it is not equivalent
to a truly primitive one.
\section{Acknowledgements}

I would like to thank the School of Mathematics of the University of Sydney for
their generosity and stimulating environment. Furthermore, I would like to
thank Don Zagier for introducing me to the subject of the paper, Bjorn Poonen
for introducing me to the method of Chabauty and Victor Flynn, Ed Schaefer and
Joe Wetherell for instructive comments and discussions.
\providecommand{\bysame}{\leavevmode\hbox to3em{\hrulefill}\thinspace}
\providecommand{\MR}{\relax\ifhmode\unskip\space\fi MR }
\providecommand{\MRhref}[2]{%
  \href{http://www.ams.org/mathscinet-getitem?mr=#1}{#2}
}
\providecommand{\href}[2]{#2}

\end{document}